\documentclass[reqno,12pt] {article}
\usepackage{amsmath, amsfonts, amssymb, amsthm,hyperref,floatrow}
\usepackage{graphicx, color,dsfont}
\usepackage{xcolor}
\usepackage{tikz}

\hypersetup{
    linktoc=page,
    linkcolor=red,          
    citecolor=blue,        
    filecolor=blue,      
    urlcolor=cyan,
    colorlinks=true           
}
\newtheorem{theorem}{Theorem}

\newtheorem{corollary}[theorem]{Corollary}

\def \E{ \mathbb E  }

\definecolor{remi}{rgb}{1,0,0}
\usepackage{titlesec}
    \titleformat{\section}[hang]
        {\color{remi}{}\bfseries\filcenter\large}
        {\thesection.}
        {0.4em}
        {}[]

\DeclareMathSymbol{\leqslant}{\mathalpha}{AMSa}{"36} 
\DeclareMathSymbol{\geqslant}{\mathalpha}{AMSa}{"3E} 
\DeclareMathSymbol{\eset}{\mathalpha}{AMSb}{"3F}     
\renewcommand{\leq}{\;\leqslant\;}                   
\renewcommand{\geq}{\;\geqslant\;}                   

\newcommand{\R}{\mathbb{R}}

\newcommand{\N}{\mathbb{N}}

\DeclareMathSymbol{\leqslant}{\mathalpha}{AMSa}{"36} 
\DeclareMathSymbol{\geqslant}{\mathalpha}{AMSa}{"3E} 
\DeclareMathSymbol{\eset}{\mathalpha}{AMSb}{"3F}     
\renewcommand{\leq}{\;\leqslant\;}                   
\renewcommand{\geq}{\;\geqslant\;}                   
\begin{document}

\title{Limiting laws of supercritical branching random walks}

   \author{Julien Barral \footnote{Universit{\'e} Paris 13, Institut Galil\'ee, LAGA, UMR CNRS 7539, 99 rue Jean-Baptiste Cl\'ement 93430 Villetaneuse, France.}, R\'emi Rhodes \footnote{Universit{\'e} Paris-Dauphine, Ceremade, UMR 7534, Place du marechal de Lattre de Tassigny, 75775 Paris Cedex 16, France.} , Vincent Vargas  \footnote{Universit{\'e} Paris-Dauphine, Ceremade, UMR 7534, Place du marechal de Lattre de Tassigny, 75775 Paris Cedex 16, France.} \footnote{The last two authors are partially supported by the CHAMU project (ANR-11-JCJC).}
}
\maketitle

\begin{abstract}
In this note, we make explicit the limit law of the renormalized supercritical branching random walk, giving credit to a conjecture formulated in \cite{BaJiRhVa} for a continuous analogue of the branching random walk. Also, in the case of a branching random walk on a homogeneous tree, we express the law of the corresponding limiting renormalized Gibbs measures, confirming, in this discrete model, conjectures formulated by physicists about the Poisson-Dirichlet nature of the jumps in the limit, and precising the conjecture by giving the spatial distribution of these jumps.\end{abstract}

\normalsize
  
\section{Introduction and notations}
 

  We consider a branching random walk on the real line $\R$. Initially, a particle sits at the
origin. Its children form the first generation; their displacements from the
origin correspond to a point process on the line. These children have children
of their own (who form the second generation), and behave-relative to their
respective positions-like independent copies of the same point process $L$. And so
on. We write $|u|=n$ if an individual $u$ is in the $n$-th generation, and denote its position $V(u)$ (and $V(e)=0$ for the initial ancestor). 
 We define the (logarithmic) moment generating function
$$\forall t\geq 0,\quad \psi(t)=\ln \mathbb{E}\Big[\sum_{|u|=1}e^{-tV(u)}\Big]\in ]-\infty;+\infty],$$
where $\mathbb E$ denotes expectation with respect to $\mathbb P$, the law of the branching random walk. 

The numbers of particles $N_n$, $n\ge 1$, in successive generations form a Galton-Watson process, which is supercritical if and only if $\ln \E(N_1)=\psi(0)>0$, a condition we assume throughout, together with  normalizing condition $\psi(1)=0$. Then the sequence
\begin{equation}
W_n=\sum_{|u|=n}e^{-V(u)},\quad n=0,1,2,\dots \big(\sum_{\emptyset}:=0\big)
\end{equation}
is a nonnegative martingale, which converges almost surely to a limit $W_\infty$. Moreover, $\mathbb{P}(W_\infty>0)>0$ if and only if $\psi'(1^-)<0$ and $\mathbb{E}\big[\big (W_1\big )\ln^+\big (W_1\big )\big ]<\infty$. In that case, the convergence holds in $L^1$ norm too and $W_\infty>0$ almost surely on the set of non-extinction of the Galton-Watson process. When $\psi'(1^-)\ge 0$, $W_\infty=0$ almost surely, and this gives rise to the issue of finding natural renormalizations of $(W_n)_{n\ge 1}$. It is straightforward to see that this reduces to considering the renormalization~of 
\begin{equation}
W_{n,\beta}=\sum_{|u|=n}e^{-\beta V(u)}\quad (\beta\ge 1),
\end{equation} 
where now $\psi'(1^-)=0$. We further introduce the derivative martingale:
\begin{equation}
Z_n=\sum_{|u|=n}V(u)e^{-  V(u)}
\end{equation}
We reinforce our assumptions. Let  $\widetilde{W}_1=W_1+\sum_{|u|=1}V(u)_+e^{- V(u)}$. We assume 
\begin{equation}\label{hyp2}
 a)\quad \psi''(1)<\infty,\quad \quad\quad  b) \quad \mathbb{E}\Big[\widetilde W_1\ln_+(\widetilde{W}_1)^3\Big] <+\infty.\end{equation}
It is known that under (\ref{hyp2}) the martingale $(Z_n)_{n\ge 1}$ converges almost surely towards a random variable $Z_\infty$, which is strictly positive conditionally on the survival of the system (see \cite{big8,aid}).
 In what follows, we investigate the limiting laws of the renormalized sequence $(W_{n,\beta})_n$ for $\beta>1$. In the following theorem, a,b) of (\ref{hyp2}) are needed for the convergence in law as it is shown in \cite{madaule}. We identify the limit under additional assumptions.
 \begin{theorem}\label{main}
Assume that  $L$ is non-lattice, $\psi(0)>0$, $\psi(1)=\psi'(1^-)=0$ and (\ref{hyp2}). Let $\beta>1$.
The random variables $(n^{\frac{3}{2}\beta}W_{n,\beta})_n$ converge in law towards a non trivial random variable $W_{\infty,\beta}$. Moreover, if there exists  $\delta>1$ such that $\mathbb{E}\Big[ \big(W_{1,\beta}+N_1\big)^{\delta}\Big] <+\infty$, the law of $W_{\infty,\beta}$ is, up to a multiplicative constant $c$, that of a stable subordinator $T_{\beta}$ of index $1/ \beta$ taken at an independent time $Z_\infty$:
$$W_{\infty,\beta}\stackrel{law}{=}cT_\beta(Z_\infty).$$
\end{theorem}
We notice that this result is obtained in \cite{Webb} in the special case where the Galton Watson-tree is homogeneous and the variables $V(u)$, $|u|=1$, are Gaussian and independent.

\section{The case of Random multiplicative cascades}\label{casc}
We point out that Theorem \ref{main} has direct applications to random multiplicative cascades since random multiplicative cascades can be seen as branching random walks where the number of children for each particle has been set to two. Random multiplicative cascades were introduced by B. Mandelbrot in \cite{mandelbrotstar} as a toy model for energy dissipation in a turbulent flow. More precisely, the closed dyadic subintervals of $[0,1]$ are naturally encoded by the nodes of the binary tree, i.e., the finite words over the alphabet $\{0,1\}$ via the mapping $u=u_1\cdots u_n\mapsto I_u=[\sum_{k=1}^nu_k2^{-k},2^{-n}+\sum_{k=1}^nu_k2^{-k}]$. 

In the notations of the previous section (recall that $\psi(1)=0$ and notice that $\psi(0)=\log(2)>0$), let us assume that $V(0)$ and $V(1)$ are i.i.d., and define a sequence of random measures on $[0,1]$ as follows. Let $dx$ denote the Lebesgue measure on $[0,1]$. We define for $n\ge 1$ the random measure
 $$\mu_n(dx)=2^n e^{-V(u(x))}\,dx,\quad \mbox{where }$$ where $u(x)$ is the unique element of $\{0,1\}^n$ such that $x\in I_{u(x)}$  when $x$ is not dyadic. We recall the following theorem of \cite{cf:KahPe} first conjectured in the seminal work \cite{mandelbrotstar}:
 
 \begin{theorem} [Kahane]
 The weak-star limit $\mu(dx)=\lim_{n\to \infty}\mu_n(dx)$ exists almost surely. Moreover, $\mu$ is almost surely non trivial if and only if $\psi'(1^-)<0$.
 \end{theorem}
  We are interested in the situation where $\psi'(1^-)\ge 0$ in which case the limiting measure vanishes. More precisely, as suggested in \cite{mandelbrotstar}, we want to find a renormalization sequence $(\lambda_n)_{n\ge 1}$ such that $\nu_n(dx)=\lambda_n \mu_n(dx)$ converges weakly, at least in law.  In accordance with the previous section we assume now that  (\ref{hyp2}) holds and $\psi'(1)=0$. We define for $\beta> 0$:
$$
\quad \mu_{n,\beta}(dx)= 2^n e^{-\beta V(u(x))}(x)\,dx.$$
When $0<\beta<1$, the situation boils down to the non-degenerate (classical) situation. For $\beta=1$, this corresponds to the so-called boundary case  which was studied in \cite{AiSh}. One can deduce from \cite{AiSh} and \cite{big8} that $n^{1/2}\mu_{n,1}$ converges in probability to a positive multiple of a positive statistically self-similar measure obtained as the almost sure weak limit of  $-\frac{d \mu_{n,\beta}}{d \beta} | _{\beta=1} (dx)$, call it $\widetilde\mu_\infty$ (the same convergence in law when $V(0)$ is Gaussian can be deduced from \cite{Webb}).


 We are especially interested in the low temperature case $\beta>1$. From Theorem \ref{main}, we deduce:  
 \begin{corollary}\label{Renorm}
Assume $\psi(1)=\psi'(1)=0$ and (\ref{hyp2}). Let $\beta>1$. Suppose there exists  $\delta>1$ such that $\psi(\delta \beta)<\infty$.   The random measures $(n^{\frac{3}{2} \beta}\mu_{n,\beta}(dx))_n$ weakly converge in law towards a random Borel measure $\mu_{\infty,\beta}(dx)$ on $[0,1]$. The law of the measure  $\mu_{\infty,\beta}(dx)$ is given by  the derivative (in the sense of distributions) of the function $t \rightarrow c T_{\beta}(\widetilde\mu_\infty[0,t])$ where $c>0$ is some positive constant and $T_{\beta}$ is a stable subordinator of index $1/ \beta$ independent from $\widetilde\mu_\infty(dx)$. 
 \end{corollary}
 
 We will denote this limiting law $c T_{\beta}(\widetilde\mu_\infty(dx))$. We stress that this implies the following result on renormalized Gibbs measures:
 \begin{corollary}\label{RenormGibbs}Under the assumptions of Corollary~\ref{Renorm}, we have the following convergence in law (for the weak convergence of measures):
$$
\lim_{n\to\infty} \frac{\mu_{n,\beta}(dx)}{\mu_{n,\beta}([0,1])} = \frac{T_{\beta}(\widetilde\mu_\infty(dx))}{T_{\beta}(\widetilde\mu_\infty[0,1])}.
$$
\end{corollary}


\section{The continuous case: Gaussian multiplicative Chaos}
In this section, we discuss the continuous analogue of the Random multiplicative cascades with lognormal weights: Gaussian multiplicative Chaos. One can work with other weights than lognormal but we stick to this case for the sake of simplicity. Gaussian multiplicative Chaos was introduced by J.P. Kahane in \cite{cf:Kah}. Gaussian multiplicative Chaos and lognormal Random multiplicative cascades bear striking similarities; in particular, one can introduce a continuous analogue to the discrete star equation (cf. \cite{AlRhVa,sohier}) which plays a key role in the identification of the limit in the above theorems. In some sense, there is a correspondance between both models in the sense that roughly  all theorems in the discrete case have continuous analogues. In the paper \cite{BaJiRhVa} (section 6 perspectives), the authors conjectured an analogue to corollary \ref{RenormGibbs} for Gaussian multiplicative Chaos. In fact, theorem \ref{main} leads to an even stronger version of the conjecture in \cite{BaJiRhVa}.

The above result on renormalized Gibbs measures, i.e. corollary  \ref{RenormGibbs}, and the conjecture in \cite{BaJiRhVa} on the limit of a continuous renormalized Gibbs measure share a common feature which is proved or conjectured to hold for many supercritical Gibbs measures (also called the glassy phase in the physics litterature) : the limiting object is a purely atomic measure whose jumps are distributed according to a Poisson-Dirichlet distribution, i.e. the ordered jumps of a renormalized stable subordinator. This universal feature (proved for instance in the case of the REM, etc.) has been conjectured by physicists (see~\cite{CarLeDou}). Corollary \ref{RenormGibbs} and the corresponding conjecture in \cite{BaJiRhVa}  are
stronger than the conjecture of physicists since it also provides the
spatial distribution of the jumps.

\section{Proofs}
\subsection{Proof of Theorem \ref{main}}
It is proved in \cite{madaule} that the sequence $(n^{\frac{3}{2}\beta}W_{n,\beta})_n$ converges in law towards a random variable $W_{\infty,\beta}$. This random variable is non trivial since
$$
n^{\frac{3}{2}\beta}W_{n,\beta}\geq \big(n^{3/2}\max_{|u|} e^{-V(u)}\big)^\beta
\geq \big(n^{3/2}e^{- \min_{|u|}V(u)}\big)^\beta
=\big(e^{- \big(\min_{|u|}V(u)-\frac{3}{2}\ln n\big)}\big)^\beta.
$$
 It is proved in \cite{aid} that the last quantity $\min_{|u|}V(u)-\frac{3}{2}\ln n$ converges in law as $n\to \infty$ towards a non trivial limit. So $W_{\infty,\beta}$ cannot be $0$. 
 It just remains to identify the limiting law $W_{\infty,\beta}$. To that purpose, it is straightforward to check that $W_{\infty,\beta}$ satisfies the star equation (\cite{mandelbrotstar}):
 $$W_{\infty,\beta}\stackrel{law}{=}\sum_{|u|=1}e^{-\beta V(u)}W_{\infty,\beta}^{(u)},$$ where conditionally on $(u,V(u),|u|=1)$ the family $(W_{\infty,\beta}^{(u)})_{|u|=1}$ are independent copies of $W_{\infty,\beta}$. This equation is carefully studied in \cite{Liu}. Since $t\mapsto \psi_\beta(t)=\ln \mathbb{E}\big[\sum_{|u|=1}e^{-\beta t V(u)}\big]$ vanishes for $t=\frac{1}{\beta}$ with $\psi_\beta'(\frac{1}{\beta})=0$, the result follows from \cite{Liu} (see also \cite{durrett} for the case of homogenous trees).\qed
 
 \subsection{Convergence of the random measures of section \ref{casc}}
It is convenient to first consider measures on the symbolic space $\{0,1\}^{\N_+}$. For each finite dyadic word $u$ we denote by $[u]$ the cylinder consisting of infinite words with $u$ as prefix and define $Q(u)=e^{-V(u)}$. $\{0,1\}^{\N_+}$ is endowed with the Borel $\sigma$-field generated by the cylinders taken as a basis of open sets. Denote the dyadic coding by $\pi: \{0,1\}^{\N_+}\to [0,1]$: $\pi (x_1x_2\cdots)=\sum_{k\ge 1}x_k2^{-k}$. Also denote by $\lambda$ the uniform measure on $\{0,1\}^{\N_+}$.

For $\beta\ge 1$ and $n\ge 1$, let $\nu_{n,\beta}$ stand for the measure whose density with respect to $\lambda$ is given by $2^nQ(u)^\beta$; we have $\mu_{n,\beta}=\nu_{\beta,n}\circ \pi^{-1}$. Also, let $\widetilde \nu_{n}$ be the signed measure whose density with respect to $\lambda$  is given by $-2^n Q(u)\log Q(u)$ over each cylinder $[u]$ of generation $n$; we have  $\widetilde \nu_{n}\circ\pi^{-1}(dx)= -\frac{d \mu_{n,\beta}}{d \beta} | _{\beta=1} (dx)$. It follows from \cite{big7} that, almost surely,  $\widetilde \nu_{n}$ is ultimately positive and converges weakly to a positive measure $\widetilde\nu_\infty$; consequently $ -\frac{d \mu_{n,\beta}}{d \beta} | _{\beta=1} (dx)$ weakly converges to $\widetilde\mu_\infty(dx)=\widetilde\nu_\infty\circ \pi^{-1}(dx)$.  Let $N^{(\beta)}_{\widetilde \nu_{\infty}}$ be a positive Borel random measure on $\{0,1\}^{\N_+}\times \R^*_+$, whose law conditionally on $\widetilde \nu_{\infty}$ is that of a Poisson point measure with intensity $\beta\, \Gamma(1/\beta)\frac{\widetilde \nu_{\infty}(dx)dz}{z^{1+1/\beta}}.$ Then define the random measure $ \widetilde\nu_{\infty,\beta} (A)= \int_A\int_{\R_+^*} z N^{(\beta)}_{\widetilde \nu_{\infty}}(dx,dz)$
and assume that $T_\beta$ is normalized so that its Laplace exponent is $-\theta^{1/\beta}$. For all continuous function $f:\{0,1\}^{\N_+}\to\R$, $|n^{3\beta/2} \nu_{n,\beta}(f)|\leq n^{3\beta/2} \|\nu_{n,\beta}\|\|f\|_\infty$ ($W_{n,\beta}=\|\nu_{n,\beta}\|$). Since the laws of the variables $n^{3\beta/2}\nu_{n,\beta}(\{0,1\}^{\N_+})$, $n\ge 1$ form a tight sequence (Theorem \ref{main}),  so do the laws of the random measures $n^{3\beta/2} \nu_{n,\beta}$, $n\ge 1$, on the compact set $\{0,1\}^{\N_+}$ for the weak convergence of measures. Let us show that the unique limit point of this sequence is the law of $ c\widetilde\nu_{\infty,\beta}$, for some $c>0$. Since the linear combinations of the indicator functions of cylinders generate the space of real-valued continuous functions on $\{0,1\}^{\N_+}$, it is enough to show that for some $c>0$, for each $p\ge 1$, $(n^{3\beta/2} \nu_{n,\beta}([u]))_{u\in\{0,1\}^p}$ converges in law to $c(\widetilde \nu_{\infty,\beta}([u]))_{u\in\{0,1\}^p}$ as $n\to\infty$. Using the self-similarity of the construction and Theorem~\ref{main}, it is not difficult to see that $(n^{3\beta/2} \nu_{n,\beta}([u]))_{u\in\{0,1\}^p}$ converges in law to $(Q(u)^\beta W_{\infty,\beta}(u))_{u\in\{0,1\}^p}$, where $(W_{\infty,\beta}(u))_{u\in\{0,1\}^p}$ is independent of $(Q(u))_{u\in\{0,1\}^p}$ and its components are independent copies of $W_{\infty,\beta}$. On the other hand, a calculation combining standard properties of Poisson measures with the fact that $(\widetilde\nu_\infty([u]))_{u\in\{0,1\}^p}=(Q(u)Z_\infty(u))_{u\in\{0,1\}^p}$,  where $(Z_{\infty}(u))_{u\in\{0,1\}^p}$ is independent of $(Q(u))_{u\in\{0,1\}^p}$ and its components are independent copies of $Z_{\infty}$, shows that the Laplace transform of $(\nu_{\infty,\beta}([u]))_{u\in\{0,1\}^p}$ equals that of $(Q(u)^\beta T_\beta^{(u)}(Z_{\infty}(u))_{u\in\{0,1\}^p}$, where the $T_\beta^{(u)}$ are independent copies of $T_\beta$, and  also are independent of the $(V(u),Z_\infty(u))_{u\in \{0,1\}^p}$. Due to Theorem~\ref{main}, this yields the desired equality in law. It follows that $n^{3\beta/2} \mu_{n,\beta}$ converges weakly in law to $c\widetilde\mu_{\infty,\beta}:= c \widetilde\nu_{\infty,\beta}\circ\pi^{-1}$. Moreover, the measure $\widetilde \nu_{\infty}$ assigns no mass to the countable set of points of $\{0,1\}^{\N_+}$ encoding dyadic points (see \cite{Ba}), hence so does $\widetilde\nu_{\infty,\beta}$, and after projection so do $\widetilde\mu_\infty$ and $\widetilde\mu_{\infty,\beta}$ with dyadic numbers. The previous information are  enough to get the equality in law, for each $p\ge 1$, of  $ (\widetilde\mu_{\infty,\beta}(I_u))_{u\in\{0,1\}^p}$ and that of $(T_\beta (\widetilde\mu_\infty)(I_u)) _{u\in\{0,1\}^p}$, and then the equality in law of $\widetilde\mu_{\infty,\beta}(dx)$ and $T_\beta (\widetilde\mu_\infty(dx))$, hence the desired conclusion.  \qed

 \hspace{1 cm}

\noindent {\bf Acknowledgements:} The authors are thankful to J.P. Bouchaud and Alberto Rosso for useful discussions and to Yueyun Hu for orienting them to Thomas Madaule's paper.

\hspace{10 cm}

\end{document}